\documentclass{elsart}
\usepackage{amssymb,amscd}
\begin{document}

\newcommand{\cC}{\mathcal{C}}
\newcommand{\cD}{\mathcal{D}}
\newcommand{\cE}{\mathcal{E}}
\newcommand{\cF}{\mathcal{F}}
\newcommand{\cG}{\mathcal{G}}
\newcommand{\cH}{\mathcal{H}}
\newcommand{\cM}{\mathcal{M}}
\newcommand{\cN}{\mathcal{N}}
\newcommand{\cV}{\mathcal{V}}
\newcommand{\ccV}{{}^c\mathcal{V}}
\newcommand{\ch}{\operatorname{ch}}
\newcommand{\coker}{\operatorname{coker}}
\newcommand{\cssX}{{}^c\!S^*X}
\newcommand{\ctX}{{}^cTX}
\newcommand{\ctsX}{{}^cT^*X}
\newcommand{\cun}{\cC^{\infty}}
\newcommand{\cuno}{\cC^{\infty}_0}
\newcommand{\dc}{\cD_c}
\newcommand{\Hom}{\operatorname{Hom}}
\newcommand{\iH}{{I_H}}
\newcommand{\iM}{{I_M}}
\newcommand{\ind}{\operatorname{index}}
\newcommand{\pc}{\Psi_{\!c}}
\newcommand{\pcs}{\Psi_{\!c,\sus}}
\newcommand{\ps}{\Psi_\sus}
\newcommand{\px}{\partial_x}
\newcommand{\re}{\operatorname{Re}}
\newcommand{\Res}{\operatorname{Res}}
\newcommand{\rc}{{\rho_c}}
\newcommand{\sD}{\sigma(D)}
\newcommand{\sign}{\operatorname{sign}}
\newcommand{\sssM}{{S_\sus^*(M)}}
\newcommand{\ssX}{S^*X}
\newcommand{\sus}{{\operatorname{sus}}}
\newcommand{\supp}{\operatorname{supp}}
\newcommand{\ta}{\tilde{a}}
\newcommand{\tD}{{\tilde{D}}}
\newcommand{\td}{{\tilde{\delta}}}
\newcommand{\tT}{{\tilde{T}}}
\newcommand{\Td}{\operatorname{Td}}
\newcommand{\Tr}{\operatorname{Tr}}
\newcommand{\tssM}{T_\sus^* M}
\newcommand{\tX}{{\tilde{X}}}

\begin{frontmatter}
\title{Cusp geometry and the cobordism invariance of the index}
\thanks{The author has been partially supported by the
Research and Training Network HPRN-CT-1999-00118 ``Geometric Analysis'' 
funded by the European Commission}
\author{Sergiu Moroianu}
\ead{moroianu@alum.mit.edu}
\ead[url]{alum.mit.edu/www/moroianu}
\date{Version of 22 June 2004}
\address{Institutul de Matematic\u{a} al Academiei Rom\^{a}ne\\
P.O. Box 1-764\\RO-70700
Bucharest, Romania}
\begin{abstract}
The cobordism invariance of the index on closed manifolds is reproved using 
the calculus $\pc$ of cusp pseudodifferential operators 
on a manifold with boundary.
More generally, on a compact manifold with corners, the existence of a symmetric cusp
differential operator of order $1$ and of Dirac type near the boundary implies
that the sum of the indices of the induced operators on the hyperfaces is null.
\end{abstract}
\end{frontmatter}

\section{Introduction}\label{intro}

Thom's discovery \cite{thom} of the cobordism invariance of the topological
signature led Hirzebruch \cite{hirz} to 
identify the signature of the intersection form of a closed oriented 
$4k$-dimensional manifold with the $L$-number constructed from the 
Pontryagin classes, in what was to become one of last century's most 
influential formulae:
\[\sign(M)=L(M).\]
Inspired by this result, Atiyah and Singer proposed in \cite{as} an 
extension of the signature formula which gave
the answer to the general index problem for elliptic operators on 
closed manifolds. Their program 
was carried out in \cite{palais}. The key ingredients of the proof were
the use of pseudodifferential operators and the cobordism invariance 
of the index of twisted signature operators. Instead of explaining what 
this is, let us state a more general result, that we will later extend 
to manifolds with corners.
\begin{thm}\label{main}
Let $M$ be a closed manifold, $\cE^\pm$ vector bundles over $M$ and 
$D:\cun(M,\cE^+)\to\cun(M,\cE^-)$
an elliptic differential operator of order $1$.
Assume that 
\begin{enumerate}
\item $M$ (not necessarily orientable) is the boundary of a compact
manifold $X$; fix a Riemannian metric on $X$ which is of product type
in a product decomposition $M\times[0,\epsilon)$ of $X$ near $M$;
\item there exists a vector bundle $\cE\to X$ such that 
$\cE|_M=\cE^+\oplus \cE^-$;  identify $\cE$ over $M\times[0,\epsilon)$
with the pull-back of $\cE^+\oplus \cE^-$ from $M$, and fix a metric on 
$\cE$ which is constant in $t\in[0,\epsilon)$ such that $\cE^+\perp \cE^-$;
\item there exists a formally self-adjoint
elliptic operator $\delta$ acting on $\cun(X,\cE)$ which near $M$ has the form
\begin{equation}\label{ea}
\delta=\left[\begin{array}{cc}-i\frac{d}{dt}&D^*\\D&i\frac{d}{dt}\end{array}\right].
\end{equation}\label{hyp}
\end{enumerate}
Then $\ind(D)=0.$
\end{thm}

Since the topological significance of $\ker D$, which played a key role 
for the signature problem, 
is lost for arbitrary elliptic operators, the 
proof from \cite[Chapter XVII]{palais} 
of even a particular case of Theorem \ref{main} had to rely
on a fairly complicated analysis of boundary value problems. 
Atiyah and Singer found later \cite{as2} a purely
$K$-theoretic proof of the index theorem, from which the cobordism 
invariance of twisted signatures follows. From a modern perspective,
Theorem \ref{main} is also a consequence of the following commutative diagram
in analytic K-homology \cite{higroe} (I am indebted to the 
referee for this remark):
\begin{equation}\label{refer}
\begin{CD}
K_{1}(X,M)@>{\partial}>>K_0(M)\\
 @VVV  @V{\ind}VV\\
K_1(\mathrm{pt},\mathrm{pt})@>>> K_0(\mathrm{pt})
\end{CD}
\end{equation}
An operator $\delta$ as in Theorem \ref{main} defines an element in 
$K_{1}(X,M)$ with $\partial(\delta)=D\in K_0(M)$. On the other hand,
$K_1(\mathrm{pt},\mathrm{pt})=0$ so $\ind(D)=0$. 
Nevertheless, there is a great deal of work in either proving the Atiyah-Singer 
formula or in constructing analytic K-homology and proving the commutativity of
(\ref{refer}).
Thus it is legitimate to ask how deep the cobordism invariance of the index
really is. Our first result is a new 
proof of Theorem \ref{main} 
by some clever manipulations with noncommutative residues
inside the calculus of cusp pseudodifferential operators on $X$
(arguably the simplest example of a pseudodifferential calculus
on a manifold with boundary,   
constructed using the theory of boundary fibration structures
of Melrose \cite{mel1}). Note that 
several proofs of Theorem \ref{main} have been obtained lately for Dirac
operators (e.g., \cite{braver}, \cite{higson}, \cite{lesch},
\cite{nicol}). 
%Our method has the advantage that it generalizes easily
%to pseudodifferential operators (Theorem \ref{pt}).
A $K$-theoretic statement of the cobordism invariance of the index  
was proved recently by Carvalho \cite{carvalho,carvart} via
the topological approach of \cite{as2}.

The main result of the paper concerns the cobordism invariance problem on 
manifolds with corners. Let $X$ be a compact manifold with corners
and $\cF_1,\ldots,\cF_k$ its boundary hyperfaces, possibly disconnected.
We refer to \cite{in} for an overview of cusp pseudodifferential operators
on manifolds with corners. 
Let $A$ be a symmetric cusp pseudodifferential operator on $X$. Under 
certain algebraic conditions which we call 
"being of Dirac type at the boundary", 
$A$ induces cusp elliptic operators $D_j$
on each hyperface $\cF_j$. We assume that these operators are fully elliptic,
which is equivalent to $D_j$ being Fredholm on suitable cusp Sobolev spaces.
Then, under the assumption that $A$ is a first-order differential operator,
we prove in Theorem \ref{maind} that 
$\sum_{j=1}^k \ind(D_j^+)=0$. 
The proof is inspired from the closed case; 
we look at a certain meromorphic function of zeta-type 
in several complex variables. A special Laurent coefficient of this function
will give on one hand the sum of the indices of $D_j$, and on the other hand 
it will vanish. 

An index formula for fully elliptic cusp operators
on manifolds with corners was given in \cite{in}. Inadvertently, we stated 
there the result only for scalar operators, however the formula applies
\emph{ad literam} to operators acting on sections of a vector bundle.
The result from Section \ref{cmain} is in a certain sense the 
odd-dimensional version of that formula. Unlike in the closed case, 
it seems difficult to obtain 
the cobordism invariance directly from the general index formula.
Note however 
that for admissible Dirac operators, Theorem \ref{maind}
can be deduced from results of Loya \cite[Theorem 8.11]{loya},
Bunke \cite[Theorem 3.14]{bunke}, and also 
from a particular case of \cite[Theorem 5.2]{in}, since 
in that case the index density is a characteristic form. We discuss this
briefly at the end of Section \ref{cmain}.

From a different point of view, Melrose 
and Rochon \cite{mr} use a slightly modified version of the cusp 
algebra to study 
the index of families of operators on manifolds with boundary. It would 
be interesting to combine their approach with ours, to treat for instance
the cobordism invariance of the families index.

In Sections \ref{rev} and \ref{pro} we will use the notation and 
some simple results from \cite{in}; we refer the reader to 
\cite{meni96c} for a thorough treatment of the cusp algebra on 
manifolds with boundary. In the second part of the paper dealing
with manifolds with corners we will rely heavily on \cite{in}. Some 
familiarity with that paper must therefore be assumed.

{\small
{\bf Acknowledgments.}
I am indebted to Catarina Carvalho for disclosing to me her 
recent result \cite{carvalho}, and to the referee for extremely helpful 
comments. I thank Maxim Braverman, John Lott, Richard Melrose and Mihai Pimsner
for useful remarks.
I also wish to thank the \emph{Equipe de G\'eom\'etrie Noncommutative de
Toulouse} for their warm hospitality at the Paul Sabatier University,
where this paper was completed.
}

\section{Review of Melrose's cusp algebra} \label{rev}
Let $X$ be a compact manifold with boundary $M$, not necessarily 
orientable, and $x:X\to\Rset_+$ a boundary-defining function 
(i.e., $M=\{x=0\}$ and $d x$ is never zero at $x=0$).
A \emph{cusp vector field} on $X$ is a smooth vector field $V$ such that 
$d x(V)\in x^2\cun(X)$ (we remind the reader that $\cun(X)$ is defined as 
the space of restrictions to $X$ of smooth functions on the double of $X$, 
or equivalently as the space of those smooth functions in the interior 
of $X$ that admit Taylor series expansions at $M$). The space of cusp vector 
fields forms a Lie subalgebra $\ccV(X)\hookrightarrow \cV(X)$, whose
universal enveloping algebra is by definition the algebra $\dc(X)$ of
scalar cusp differential operators. Moreover $\ccV(X)$ is a finitely 
generated projective $\cun(X)$-module (in a product decomposition 
$M\times[0,\epsilon)\hookrightarrow X$, a local basis is given by
$\{x^2\frac{\partial}{\partial x}, \frac{\partial}{\partial y_j}\}$ 
where $y_j$ are local coordinates on $M$).
Thus by the Serre-Swan theorem there exists a vector bundle $\ctX\to X$ 
such that $\ccV(X)=\cun(X,\ctX)$. 

A cusp differential operator of positive order can never be elliptic at $x=0$. 
Nevertheless, there exists a natural cusp principal symbol map  
surjecting onto the smooth polynomial functions on $\ctsX$ of 
homogeneity $k$, $\sigma:\dc^k(X)\to \cun_{[k]}(\ctsX)$. 
A cusp operator will therefore be called elliptic if 
its principal symbol is invertible on $\ctsX\setminus\{0\}$. 

For any vector bundles $\cF,\cG$ over $X$ let 
\[\dc(X,\cF,\cG):=\dc(X)\otimes_{\cun(X)} \Hom(\cF,\cG).\]
It is straightforward to extend the definition of $\sigma$ to the bundle case.

\subsection{Example}\label{exam} Assume that the hypothesis of Theorem 
\ref{main} is fulfilled. The metric on $X$ is a product metric near $M$,
\[g^X=dt^2+g^M.\]
Extend $\delta$ to the manifold $\tX=M\times(-\infty,0)\cup X$, obtained by 
attaching an infinite cylinder to $X$,
by Eq.\ (\ref{ea}). Let $\psi:X^\circ\to\tX$
be any diffeomorphism extending
\[M\times(0,\epsilon)\ni(y,x)\mapsto(y,-\frac1x)\in\tX.\]
Then the pull-back of $\delta$ through $\psi$ takes the form
\[A:=\psi^*\delta=\left[\begin{array}{cc}-ix^2\frac{d}{dx}&D^*\\D&
ix^2\frac{d}{dx}\end{array}\right]\]
since $t=-\frac1x$ near $x=0$. Thus, $A$ is a cusp differential operator. Moreover,
$A$ is symmetric with respect to the cusp metric $\psi^* g^X$,
which near $M$ takes the form
\[g_c^X=\frac{dx^2}{x^4}+g^M.\]
The metric $g_c^X$ is degenerate at $x=0$, however it is non-degenerate 
as a cusp metric in the sense that it induces 
a Riemannian metric on the bundle $\ctX\to X$. The operator $A$ is elliptic 
(in the cusp sense) and acts as an unbounded operator in $L^2_c(X,\cE)$, the space 
of square-integrable sections in $\cE\to X$ with respect to the metric $g_c^X$.

\subsection{The indicial family} This is a "boundary symbol" map, associating 
to any cusp operator $P\in\dc(X,\cE,\cF)$ a family of differential
operators on $M$ with one real polynomial parameter $\xi$ as follows:
\[\iM(P)(\xi)=\left(e^{\frac{i\xi}{x}}P e^{-\frac{i\xi}{x}}\right)|_M\]
where restriction to $M$ is justified by the mapping properties
\begin{eqnarray*}
P:\cun(X,\cE)&\to&\cun(X,\cF)\\
P:x\cun(X,\cE)&\to&x\cun(X,\cF)
\end{eqnarray*}
and by the isomorphism $\cun(M)=\cun(X)/x\cun(X)$.

From the definition we see directly for the cusp operator $A$ 
constructed in Example \ref{exam} that
\begin{equation}\label{cna}
\iM(A)=\left[\begin{array}{cc}\xi&D^*\\D&-\xi\end{array}\right].
\end{equation}

Ellipticity does not make $A$ Fredholm on $L^2_c(X,\cE)$, 
essentially because the Rellich lemma does not hold on non-compact domains. 
To apply a weighted form of the Rellich lemma we need an extra property, 
the invertibility of the indicial family 
$\iM(A)$ for all values of the parameter $\xi$;
thus $A$ is Fredholm precisely when $D$ is invertible, see 
\cite[Theorem 3.3]{in}. Elliptic cusp 
operators with invertible indicial family are called fully elliptic.

\subsection{Cusp pseudodifferential operators} By a micro-localization process
one constructs \cite{meni96c} a calculus of pseudodifferential operators
$\pc^\lambda(X)$, $\lambda\in\Cset$, in which $\dc(X)$ sits as the 
subalgebra of differential operators (the symbols used in the definition 
are classical of order $\lambda)$. By composing with the multiplication
operators $x^z$, $z\in\Cset$, we get a calculus with two indices
\[\pc^{\lambda,z}(X,\cF,\cG):=x^{-z}\pc^\lambda(X,\cF,\cG)\]
such that $\pc^{\lambda,z}(X,\cE,\cF)\subset \pc^{\lambda',z'}(X,\cE,\cF)$ 
if and only if $\lambda'-\lambda\in\Nset$ and $z'-z\in\Nset$ (since we work 
with classical symbols). Also,
\[\pc^{\lambda,z}(X,\cG,\cH)\circ\pc^{\lambda',z'}(X,\cF,\cG)\subset 
\pc^{\lambda+\lambda',z+z'}(X,\cF,\cH).\]

By closure, cusp operators act 
on a scale of weighted Sobolev spaces $x^\alpha H_c^\beta$:
\[\pc^{\lambda,z}(X,\cF,\cG)\times x^\alpha H_c^\beta(X,\cF)\to x^{\alpha-\Re(z)}
H_c^{\beta-\Re(\lambda)}(X,\cG).\]
The principal symbol map and the indicial family extend to multiplicative
maps on $\pc(X)$. The indicial family takes values in the space $\Psi_\sus(M)$ 
of families of operators on $M$ with one real parameter $\xi$, with joint 
symbolic behavior in $\xi$ and in the cotangent variables of $T^*M$
($1$-suspended pseudodifferential operators in the terminology 
of \cite{meleta}).

The following result gives a hint of what families of operators actually define
suspended operators.

\begin{lem} \label{prod}
Let $z,w\in\Cset\cup\{-\infty\}$, $P\in\Psi^{z}(M)$ and $\phi\in\cun(\Rset)$.
Then $\xi\mapsto \phi(\xi)P$ belongs to $\ps^w(M)$ if and only if one of 
the following two conditions is fulfilled:
\begin{enumerate}
\item $z=-\infty$ and $\phi$ is a rapidly decreasing 
(i.e., Schwartz) function.
\item $P$ is a differential operator and $\phi$ is a polynomial.
\end{enumerate}
In the first case $w=-\infty$, while in the second case $w=z+\deg(\phi)$.
\end{lem}

\subsection{Analytic families of cusp operators}\label{afoco}
Let $Q\in\pc^{1,0}(X,\cE)$ be a positive fully elliptic cusp operator of order 
$1$. Then the complex powers $Q^\lambda$ form an analytic family of 
cusp operators of order $\lambda$ (this is proved using Bucicovschi's 
method \cite{buc}). 

Let $\Cset^2\ni(\lambda,z)\mapsto P(\lambda,z)\in\pc^{\lambda,z}(X,\cE)$ 
be an analytic family in two complex variables. Then $P(\lambda,z)$
is of trace class (on $L^2_c(M,\cE)$) for 
\[\Re(\lambda)<-\dim(X),\quad \Re(z)<-1,\]
and $\Tr(P)$ is analytic there as a function of $(\lambda,z)$. Moreover,
$\Tr(P)$ extends to $\Cset^2$ meromorphically with at most simple poles 
in each variable at $\lambda\in\Nset-\dim(X)$, $z\in\Nset-1$. By analogy with 
the Wodzicki residue, we can give a formula for the residue at $z=-1$ as a 
meromorphic function of $\lambda$ (this is essentially 
\cite[Proposition 4.5]{in}).
\begin{prop}\label{p54}
Let $\Cset^2\ni(\lambda,z)\mapsto P(\lambda,z)\in\pc^{\lambda,z}(X,\cE)$ 
be an analytic family. Then 
\begin{equation}\label{trz}
\Res_{z=-1}\Tr(P(\lambda,z))=\frac{1}{2\pi}\int_\Rset \Tr(\iM(x^{-1}P(\lambda,-1)))
d\xi.
\end{equation}
\end{prop}
\begin{pf} 
The trace on the right-hand side is the trace on $L^2(M,\cE_{|M})$. 
Both terms are meromorphic functions in $\lambda\in\Cset$.
By unique continuation, it is thus enough to prove the identity for 
$\Re(\lambda)<-\dim(X)$. We write both traces as the integrals of the
trace densities of the corresponding operators, i.e., the restriction of their
distributional kernel to the diagonal. 

For any vector bundle $V\to X$ we denote by $\Omega(V)\to X$ the associated 
density bundle.
Let $F:\Cset\to \cun(X,\Omega(\ctX))$ be a holomorphic family of smooth 
cusp-densities. Then $x^2F(z)\in\cun(X,\Omega(TX))$, and hence
$z\mapsto \int_X x^{-z}F(z)$ is holomorphic for $\Re(z)<-1$; moreover, 
its residue at $z=-1$ is easily seen to equal 
$\int_M (x^2\px\lrcorner F(-1))_{|M}$.
We apply this fact to
the trace density of $P(\lambda,z)$ multiplied with $x^z$. We view $M$
as the intersection of the cusp diagonal with the cusp front face inside
the cusp double space $X^2_c$ (see \cite{in}). Recall from \cite{in} 
or \cite{meni96c} that the cusp front face is the total space of a real line
bundle, and $M$ lives inside the zero section. The indicial family is obtained
by restricting a Schwartz kernel to the front face, then Fourier transforming 
along the fibers. The result follows from the Fourier inversion formula
\[f(0)=\frac{1}{2\pi}\int_\Rset \hat{f}(\xi)d\xi\]
applied in the fibers of the cusp front face over $M$.
\end{pf}

\section{Cobordism invariance on manifolds with boundary}\label{pro}

The self-contained proof of the cobordism invariance of 
the index given below serves as a model for the general 
statement on manifolds with corners. 
\begin{pf*}{Proof of Theorem \ref{main}.}
We have seen in Example \ref{exam} that the hypothesis of Theorem 
\ref{main} is equivalent to the existence of an 
elliptic symmetric cusp operator $A$ satisfying (\ref{cna}).

Let $\phi:\Rset\to\Rset$ be a non-negative Schwartz function with $\phi(0)=1$.
Let $P_{\ker D}\in\Psi^{-\infty}(M,\cE^+)$, 
$P_{\coker D}\in\Psi^{-\infty}(M,\cE^-)$ be the (finite-rank) 
orthogonal projections on the kernel and cokernel of $D$. These projections 
belong to $\Psi^{-\infty}(M)$ by elliptic regularity.
By Lemma \ref{prod},
\[r(\xi):=\left[\begin{array}{cc}\phi(\xi)P_{\ker D}&0\\0&\phi(\xi)P_{\coker D}
\end{array}\right]\] 
belongs to $\ps^{-\infty}(M,\cE^+\oplus\cE^-)$ and is non-negative, so 
it is the indicial family of a non-negative
cusp operator $R\in\pc^{-\infty,0}(X,\cE)$.
By (\ref{cna}), 
\begin{equation}\label{cnar}
\iM(A^2+R)=\left[\begin{array}{cc}\xi^2+D^*D+\phi(\xi)P_{\ker D}&0\\0&
\xi^2+DD^*+\phi(\xi)P_{\coker D}\end{array}\right]
\end{equation}
so $A^2+R$ is fully elliptic. It follows that 
$P_{\ker(A^2+R)}\in\pc^{-\infty,-\infty}(X,\cE)$ (by elliptic regularity 
with respect to the two symbol structures) so $A^2+R+P_{\ker(A^2+R)}$
is a positive cusp operator. Finally set 
\[Q:=(A^2+R+P_{\ker(A^2+R)})^{1/2}\]
and let $Q^\lambda$ be its complex powers. 
Note that $Q^2-A^2\in \pc^{-\infty,0}(X,\cE)$ so 
\begin{equation}\label{ci}
[A,Q^\lambda]\in\pc^{-\infty,0}(X,\cE).
\end{equation}

Let $P(\lambda,z)\in\pc^{-\lambda -1, -z-1}(X,\cE)$ be the analytic family 
of cusp operators 
\[P(\lambda,z):=[x^z,A]Q^{-\lambda -1}.\]
From the discussion in Subsection \ref{afoco}, $\Tr(P(\lambda,z))$
is holomorphic in $\{(\lambda,z)\in\Cset^2; \Re(\lambda>\dim(X)-1, \Re(z)>0\}$
and extends meromorphically to $\Cset^2$. We keep the notation $\Tr(P(\lambda,z))$
for this extension. Note that although $P(\lambda,0)=0$, there is no reason 
to expect the meromorphic extension
$\Tr(P(\lambda,z))$ to vanish at $z=0$; rather, $\Tr(P(\lambda,z))$
will be regular in $z$ near $z=0$.
Our proof of Theorem \ref{main} will consist of computing in two 
different ways the complex number 
\[N:=\Res_{\lambda=0} \left(\Tr(P(\lambda,z))|_{z=0}\right)\]
where $(\cdot)|_{z=0}$ denotes the regularized value in $z$ 
at $z=0$, which is a meromorphic function of $\lambda$. In other words, 
$N$ is the coefficient of $\lambda^{-1}z^0$ in the Laurent expansion 
of $\Tr(P(\lambda,z))$ around $(0,0)$. Evidently, we can also 
take the residue in $\lambda$ \emph{before} evaluating at $z=0$; 
in that case, the output of the residue is a meromorphic function in $z$.

On one hand, we claim that
\[ \Tr(P(\lambda,z))=\Tr(x^z[A,Q^{-\lambda -1}]) \]
for all $\lambda,z\in\Cset$. Since $[x^z,A]Q^{-\lambda -1}
=x^z[A,Q^{-\lambda -1}]+[x^zQ^{-\lambda-1},A]$, the claim is equivalent 
to showing that the meromorphic function 
\[(z,\lambda)\mapsto \Tr([x^zQ^{-\lambda-1},A])\]
vanishes identically. Indeed, for $\Re(\lambda)>\dim(X), \Re(z)>1$ 
this vanishing holds by the trace property, and 
unique continuation proves the claim in general.
Furthermore, $x^z[A,Q^{-\lambda -1}]\in\pc^{-\infty,-z}(X,\cE)$ by (\ref{ci})
so in fact 
\[(\lambda,z)\mapsto\Tr(x^z[A,Q^{-\lambda -1}])=\Tr(P(\lambda,z))\] 
is analytic in $\lambda\in\Cset$. In conclusion 
$\Tr(P(\lambda,z))$ is regular in 
$\lambda$ at $\lambda=0$, so 
\begin{equation}\label{no}
N=0.
\end{equation}

On the other hand, $P(\lambda,0)=0$ so 
\[U(\lambda,z):=z^{-1}P(\lambda, z)\in 
\pc^{-\lambda -1, -z-1}(X,\cE)\]
is an analytic family. Set
$[\log x,A]:=(z^{-1}[x^z,A])|_{z=0}\in\pc^{0,1}(X,\cE)$. Then 
$U(\lambda, 0)=[\log x,A]Q^{-\lambda-1}$ and 
\begin{eqnarray*}
\iM(x^{-1}U(\lambda,0))&=&\iM(x^{-1}[\log x,A])\iM(Q^{-\lambda-1})\\
&=&\left[\begin{array}{cc}i&0\\0&-i\end{array}\right]
\iM(A^2+R+P_{\ker(A^2+R)})^{-\frac{\lambda+1}{2}}
\end{eqnarray*}
where $\iM(A^2+R+P_{\ker(A^2+R)})$ is given by (\ref{cnar}) because 
$\iM(P_{\ker(A^2+R)})=0$. Using (\ref{trz}) we get
\begin{eqnarray}
\Tr(P(\lambda, z))|_{z=0}&=&\Res_{z=0}\Tr(z^{-1}P(\lambda, z))\nonumber\\
&=&\frac{1}{2\pi}\int_\Rset \Tr(\iM(x^{-1}(U(\lambda,0)))d\xi\nonumber\\
&=&\frac{i}{2\pi}\int_\Rset 
\left(\Tr(\xi^2+D^*D+\phi(\xi)P_{\ker D})^{-\frac{\lambda+1}{2}}\right.
\nonumber\\
&&\left. -\Tr(\xi^2+DD^*+\phi(\xi)P_{\coker D})^{-\frac{\lambda+1}{2}}
\right)d\xi\label{nd}
\end{eqnarray}
For each fixed $\xi$ we compute the trace using an orthonormal 
basis of $L^2(X,\cE^\pm)$ given
by eigensections of $D^*D$, respectively $DD^*$. Clearly the 
contributions of nonzero eigenvalues cancel in (\ref{nd}) so we are left with
\[\Tr(P(\lambda, z))|_{z=0}=\ind(D)\int_\Rset
(\xi^2+\phi(\xi))^{-\frac{\lambda+1}{2}}d\xi.\]
The reader will easily convince herself that the residue 
\[\Res_{\lambda=0}\int_\Rset (\xi^2+\phi(\xi))^{-\frac{\lambda+1}{2}}d\xi\]
is independent of the Schwartz function $\phi$ and equals $2$. Thus
\begin{eqnarray*}
N&=&\Res_{\lambda=0}\Tr(P(\lambda, z))|_{z=0}\\
&=&\frac{i}{\pi}\ind(D).
\end{eqnarray*}
Together with (\ref{no}) this finishes the proof of Theorem \ref{main}.
\end{pf*}

\section{Cobordism invariance on manifolds with corners} \label{cmain}

Let $X$ be a manifold with corners in the sense
of Melrose \cite{mel1}. Let $\cM_1(X)$ be the set of boundary hyperfaces, 
possibly disconnected,
and for $H\in\cM_1(X)$ fix $x_H$ a defining function for $H$.
We fix a product cusp metric $g^X$ on the interior of $X$, 
which means iteratively that near each hyperface $H$, $g^X$ 
takes the form 
\[g^X=\frac{dx_H^2}{x_H^4}+g^H\]
for a product cusp metric on $H$. The algebra of cusp differential operators
on $X$ is simply the universal enveloping algebra of the Lie algebra
of smooth vector fields on $X$ of finite length with respect to $g^X$.
The algebra $\pc(X)$ of cusp pseudodifferential operators
was described in \cite{in} 
(see for instance \cite{mel1} for the general ideas behind such constructions).
In this section the reader is assumed to be familiar with \cite{in}.
Our main result is inspired from Theorem \ref{main}. 

\begin{thm} \label{maind}
Let $X$ be a compact manifold with corners and 
\[D_H:H^1_c(H,\cE_H^+)\to L^2_c(H,\cE_H^-)\] 
a fully elliptic cusp differential operator of order $1$ for each 
hyperface $H$ of $X$.
Assume that there exists a hermitian vector bundle $\cE\to X$ 
with product metric near the corners and 
$A\in\pc^1(X,\cE)$ a (cusp) elliptic symmetric 
differential operator, such that for each $H\in\cM_1(X)$,
$\cE|_H\cong \cE_H^+\oplus\cE_H^-$ and 
\begin{equation}\label{dtb}
\iH(A)(\xi_H)=
\left[\begin{array}{cc}\xi_H&D_H^*\\D_H&-\xi_H\end{array}\right].
\end{equation}
Then 
\[\sum_{H\in\cM_1(X)}\ind(D_H)=0.\]
\end{thm}
\begin{rem}\label{rems}
The existence of $A$ requires the following compatibility condition for 
$D_H$, $D_G$ near $H\cap G$: 
\begin{eqnarray*} I_G(D_H)&=&i\xi_G+D_{HG}\\I_H(D_G)&=&i\xi_H-D_{HG}
\end{eqnarray*}
where $D_{HG}$ is a symmetric invertible differential operator on $\cE_H^+$
over the corner $H\cap G$, and $\cE_H^+$, $\cE_H^-$, $\cE_G^+$, 
$\cE_G^-$ are identified over $G\cap H$ by elementary linear algebra. 
We say that $A$ satisfying (\ref{dtb}) is \emph{of Dirac type near 
the boundary}, since the spin Dirac operator on a manifold with 
corners satisfies this condition.
\end{rem}
\begin{pf*}{Proof of Theorem \ref{maind}} For each $H\in\cM_1(X)$, the operator $I_H(A)$ is fully elliptic 
(as a suspended cusp operator), however it is invertible if and only if 
$D_H$ is invertible (this is seen easily by looking at the diagonal operator
$I_H(A)^2$). We are interested exactly in the case when $D_H$ 
has non-zero index, thus typically $A$ is not fully elliptic. Nevertheless,
$D_H$ is Fredholm and its kernel is made of smooth sections vanishing rapidly
to the boundary faces of $H$. Equivalently, the orthogonal projection
$P_{\ker D_H}$ belongs to $\pc^{-\infty,-\infty}(H,\cE_H^+\oplus \cE_H^-)$.

Let $\phi$ a cut-off function with support 
in $[-\epsilon,\epsilon]$, $\phi\geq 0$, $\phi(0)=1$. Then 
(see Lemma \ref{prod})
\[ r_H(\xi_H)=\left[\begin{array}{cc}\phi(\xi_H)P_{\ker D_H}\\
&\phi(\xi_H)P_{\coker D_H}\end{array}\right] \]
belongs to $\pcs^{-\infty,-\infty}(H,\cE^+\oplus\cE^-)$ and is non-negative.
Clearly $\iH(A)^2+r_H$ is invertible; let $R_H\in\pc^{-\infty,0}(X)$
with $\iH(R_H)=r_H$, $I_G(R_H)=0$ for $G\neq H$ (possible since $I_G(r_H)=0$)
and $R_H^*=R_H\geq 0$. Let $R:=\sum_{H\in\cM_1(X)}R_H$.
Then
$A^2+R\geq 0$ is fully elliptic; by elliptic regularity,
$P_{\ker(A^2+R)}$ belongs to $\pc^{-\infty,-\infty}(X)$.
Finally, we set 
\[Q:=(A^2+R
+P_{\ker(A^2+R)})^{1/2}.\]
The crucial property of the invertible operator $Q$ is that its complex 
powers, like $Q$ itself, commute with $A$ modulo $\pc^{-\infty,0}(X,\cE)$.

Look at the function
\[\Cset^{\cM_1(X)}\times\Cset\ni(z,\lambda)\mapsto N(z,\lambda):=
\Tr(x^z[A,Q^{-\lambda-1}]),\]
where $x^z:=x_{H_1}^{z_{H_1}}\cdot\ldots\cdot x_{H_k}^{z_{H_k}}$. Here 
we have fixed an order on the set $\cM_1(X)=\{H_1,\ldots,H_k\}$.
By a general argument \cite[Proposition 4.3]{in}, such a function
can have at most simple poles in each of the complex variables, 
occurring at certain integers. But
the family of operators involved is of order $-\infty$ with respect
to the operator order (because of the commutativity modulo  
$\pc^{-\infty,0}(X,\cE)$). Thus in fact there is no pole in $\lambda$
at $\lambda=0$, in particular
\[\Res_{\lambda=0} N(z,\lambda)_{z=0}=0.\]
For $\Re(z_H)>-1,\Re(\lambda)>\dim(X)$ the trace property allows us to write
\begin{eqnarray*}
N(z,\lambda)&=&\sum_{j=1}^k \Tr(x_{H_1}^{z_{H_1}}\cdot\ldots
[A,x_{H_j}^{z_{H_j}}]\ldots \cdot x_{H_k}^{z_{H_k}} Q^{-\lambda-1})\\
&=:&\sum_{H\in\cM_1(X)}N_H(z,\lambda).
\end{eqnarray*}
By unique continuation this identity holds for all $z,\lambda$. 
Each term $N_H(z,\lambda)$ 
of the right-hand side is a meromorphic function with at most
simple poles in each variable. In fact, $N_H(z,\lambda)$ is regular in $z_H$
at $z_H=0$, since $[A,x_H^{z_H}]$ vanishes
when $z_H=0$. By \cite[Proposition 4.5]{in} 
(see also Proposition \ref{p54}), $N_H(z,\lambda)_{|z_H=0}$
is given by 
\begin{equation}\label{fv}
\frac{1}{2\pi}\int_\Rset \Tr(x_{H_1}^{z_{H_1}}\cdot \ldots
\iH(x_H^{-1}[\log x_{H_1},A])\ldots\cdot x_{H_k}^{z_{H_k}} \iH(Q)^{-\lambda-1})
d\xi_H.
\end{equation}
By  \cite[Lemma 3.4]{in}, 
\[\iH(x_H^{-1}[\log x_H,A])=\frac1i 
\frac{\partial \iH(A)(\xi_H)}{\partial\xi_H}.\]
Now $\iH(Q)$ is a diagonal matrix, so the trace from formula (\ref{fv}) 
can be decomposed using the splitting of $\cE_{|H}$. For the terms coming from 
$\cE^\pm$, notice that the corresponding coefficient in 
$\frac{\partial \iH(A)(\xi_H)}{\partial\xi_H}$ has the 
pleasant property of being central, since it equals $\pm i$. 
Let
\begin{eqnarray*}
\widehat{x_H}^{z}&:=&x^z/x_H^{z_H},\\
T_H^+(\xi_H)&:=&(D_H^*D_H+\xi_H^2+\phi(\xi)P_{\ker D_H})^{-\frac{1}{2}},\\
T_H^-(\xi_H)&:=&(D_HD_H^*+\xi_H^2+\phi(\xi)P_{\ker D_H^*})^{-\frac{1}{2}}.
\end{eqnarray*}
With this notation we get
\begin{equation}\label{nhzz}
N_H(z,\lambda)_{|z_H=0}=\frac{i}{2\pi}\int_\Rset \Tr(\widehat{x_H}^{z}(
T_H^+(\xi_H)^{-\lambda-1}-T_H^-(\xi_H)^{-\lambda-1}))d\xi_H.
\end{equation}
The trace functional and $\widehat{x_H}^{z}$ are independent of $\xi_H$, 
so we commute them out of the integral.
We use now the identity
\[D_HT_H^+(\xi_H)^w =T_H^-(\xi_H)^wD_H\]
valid for every $w\in\Cset$, to decompose
\[T_H^-(\xi_H)^{-\lambda-1} =D_HT_H^+(\xi_H)^{-\lambda-1} 
T_H^+(0)^{-2}D_H^*+T_H^+(\xi_H)^{-\lambda-1}P_{\ker D_H^*}\]
in its components acting on $(\ker D_H^*)^\perp$, $\ker D_H^*$.
Thus 
\begin{eqnarray}
\int_\Rset T_H^-(\xi_H)^{-\lambda-1}d\xi_H&=&
D_H T_H(0)^{-\lambda-2} D_H^*\int_\Rset(1+\xi^2)^{-\lambda-1}d\xi\nonumber\\
&&+P_{\ker D_H^*} \int_\Rset (1+\phi(\xi_H))^{-\frac{\lambda+1}{2}}d\xi_H.
\label{au}
\end{eqnarray}
Similarly 
\[T_H^+(\xi_H)^{-\lambda-1} =T_H^+(\xi_H)^{-\lambda-1}
T_H^+(0)^{-2}D_H^*D_H+T_H^+(\xi_H)^{-\lambda-1}P_{\ker D_H^*}\]
so 
\begin{eqnarray}
\int_\Rset T_H^+(\xi_H)^{-\lambda-1}d\xi_H&=&
T_H(0)^{-\lambda-2} D_H^*D_H\int_\Rset(1+\xi^2)^{-\lambda-1}d\xi\nonumber\\
&& +P_{\ker D_H}\int_\Rset (1+\phi(\xi_H))^{-\frac{\lambda+1}{2}}d\xi_H.
\label{ad}
\end{eqnarray}

We are interested in the residue $\Res_{\lambda=0}N_H(z,\lambda)|_{z=0}$.
Using (\ref{au}), (\ref{ad}) we isolate in (\ref{nhzz})
the contribution of the projectors on the kernel and cokernel of $D_H$, 
and then evaluate at $\widehat{z_H}=0$. Note that these projectors belong  
to the ideal $\pc^{-\infty,-\infty}(H)$ so their contribution is regular in
$\widehat{z_H}$. Now the trace of a projector equals the dimension of its image,
while the residue at $\lambda=0$ of $\int_\Rset (1+\phi(\xi_H))^{-\frac{\lambda+1}{2}}d\xi_H$
has been seen to be $2$. Hence the contribution of the projector terms
equals $\frac{i}{\pi}\ind(D_H)$.

The function $\lambda\mapsto\int_\Rset(1+\xi^2)^{-\lambda-1}d\xi$
is regular at $\lambda=0$ with value $\pi$.
We still need to examine $\sum_{H\in\cM_1(X)} \Res_{\lambda=0}L_H(0,\lambda)$, 
where
\[L_H(\widehat{z_H},\lambda):=\Tr(\widehat{x_H}^{z}
[(D_H^*D_H+P_{\ker D_H})^{-\frac{\lambda}{2}-1}D_H^*,D_H]).\] 
The residue in $\lambda$ of $L_H$ at $\widehat{z_H}=0$ does not vanish directly, 
as one might hope at this point. We write as before 
(using the trace property for large real parts and then invoking 
unique continuation)
\begin{eqnarray*}
L_H(\widehat{z_H},\lambda)&=&\sum_{l=1}^{k-1} \Tr\left(x_{G_1}^{z_{G_1}}\ldots
[D_H,x_{G_l}^{z_{G_l}}]\ldots x_{G_{k-1}}^{z_{G_{k-1}}}\right.\\&&
\left.\cdot(D_H^*D_H+P_{\ker D_H})^{-\frac{\lambda}{2}-1}D_H^*\right)\\
&=:&\sum_{G\in\cM_1(X)\setminus \{H\}}L_{HG}(\widehat{z_H},\lambda).
\end{eqnarray*}
We see that $L_{HG}$ is regular in $z_G$ at $z_G=0$ since it is the trace of 
an analytic family of operators which vanishes at $z_G=0$. Moreover, we can
write down the value $L_{HG}(\widehat{z_H},\lambda)|_{z_G=0}$ using Proposition
\ref{p54} (or rather \cite[Proposition 4.5]{in}, its analog for 
higher codimensions). By \cite[Lemma 3.4]{in} and from Remark \ref{rems}, 
\[I_G(x_G^{-1}
[D_H, \log x_G])=i\frac{\partial I_G(D_H)(\xi_G)}{\partial \xi_G}=-1.\]
Now $I_G(P_{\ker D_H})=0$ (by full ellipticity of $D_H$) while
$I_G(D_H)(\xi_G)=i\xi_G+D_{HG}$. The term $i\xi_G$ contributes 
an odd integral in $\xi_G$ to $L_{HG}|_{z_G=0}$, so 
\[L_{HG}|_{z_G=0}=-\frac{1}{2\pi}\int_\Rset \Tr(\widehat{x_{HG}}^{z_{HG}}
(D_{HG}^2+\xi_G^2)^{-\frac{\lambda}{2}-1}D_{HG}) d\xi_G.\]
The argument is finished by Remark \ref{rems}: indeed, modulo a conjugation,
$D_{HG}=-D_{GH}$ so
the above integrand appears again with opposite sign in
$L_{GH}|_{z_H=0}$, once we replace 
the variables of integration $\xi_G, \xi_H$ with a more neutral $\xi$. 
In other words, 
$\sum_{H,G\in\cM_1(X)}L_{HG}(\widehat{z_H},\lambda)|_{\widehat{z_H}=0}=0$,
which together with our discussion on the projectors on the kernel of 
$D_H$ shows that 
\[\sum_{H\in\cM_1(X)} \Res_{\lambda=0}
N_H(z,\lambda)|_{z=0} =\sum_{H\in\cM_1(X)}\ind(D_H).\]
The left-hand side is just  $\Res_{\lambda=0}N(z,\lambda)|_{z=0}$, which 
was seen to vanish.
\end{pf*}

Recall from \cite{in} that the index of $D_H$ can be written as the (regularized)
integral on $H$ of a density depending on the full symbol of $D_H$, plus 
contributions from each corner of $H$. In the case of differential 
operators of order $1$ only the hyperfaces of $H$ have a 
non-zero contribution, which is to be 
thought of as some sort of eta invariant. The eta invariant is sensitive to 
the orientation; in our case this means that $G\cap H$ contributes to the 
index of $D_H$ and $D_G$ the same quantity with opposite signs. Thus only
the local index density detects whether our family of operators $\{D_H\}$
is cobordant to $0$ or not. If we work with twisted Dirac operators, 
the local index density is given by a characteristic form with compact support 
away from the boundary of $H$. The existence of $A$ as in Theorem \ref{maind}
ensures that this characteristic form is the restriction of a 
characteristic form from $X$ to the hyperfaces. Thus we deduce Theorem 
\ref{maind} in this case by Stokes formula. An index formula on
a manifold with corners $H$ was given in \cite{bunke} for $b$-Dirac operators, 
in \cite{loya}
for $b$-differential operators of order $1$ and in \cite{in} for cusp 
pseudodifferential operators. In this last paper the formula as stated 
covers scalar operators, but in reality it applies to operators acting on the 
sections of a vector bundle over $H$. This includes the case of Dirac operators
if, surprisingly, the boundary of $H$ is not empty. Indeed, in that case 
there exists a non-zero vector field on $H$ which identifies, 
via the principal symbol map, any two bundles (e.g., the positive and 
negative spinor bundles) related by an elliptic operator.

\section{A conjecture}

We conjecture that Theorem \ref{maind} remains true for cusp pseudodifferential 
operators of order $1$ of 
Dirac type near the boundary, in the sense of Remark \ref{rems}.
In this generality, our proof breaks down for instance
when integrating with respect to $\xi_H$.
For differential operators, we managed to show that the errors are 
concentrated at codimension $2$ corners, and cancel each other.
This seems not possible to do in the general case.
One way to proceed would be to consider cusp 
operators of order $(1,1)$, obtained by multiplying $A$ with the inverse
of the boundary defining functions. Then a power of such an operator 
of sufficiently small real part would be of trace class. Unfortunately, 
other complications arise, for instance the meromorphic extension of 
such a trace will have poles of higher order.

\bibliographystyle{plain}

\end{document}